\documentclass[preprint,12pt]{elsarticle}

\usepackage{amssymb}
\usepackage{amsthm,url}
\usepackage{amsmath}

\newtheorem{thm}{Theorem}[section]
\newtheorem{lem}[thm]{Lemma}

\newtheorem{rmk}{Remark}[section]
\numberwithin{equation}{section}

\newproof{pf}{Proof}

\journal{    }

\begin{document}

\begin{frontmatter}



\title{Nonlinear Convection in Reaction-Diffusion Equations under dynamical boundary conditions}

%
 \author{Ga\"elle Pincet Mailly \& Jean-Fran\c cois Rault }


\begin{abstract}
We investigate blow-up phenomena for positive solutions of nonlinear reaction-diffusion equations including a nonlinear convection term $\partial_t u = \Delta u - g(u) \cdot \nabla u + f(u)$ in a bounded domain of $\mathbb{R}^N$ under the dissipative dynamical boundary conditions $\sigma \partial_t u + \partial_\nu u =0$. Some conditions on $g$ and $f$ are discussed to state if the positive solutions blow up in finite time or not. Moreover, for certain classes of nonlinearities, an upper-bound for the blow-up time can be derived and the blow-up rate can be determinated.
\end{abstract}


\end{frontmatter}


\section{Introduction}
We consider the following nonlinear parabolic problem
\begin{eqnarray}\label{CRD}
\left\{
    \begin{array}{ll}
			\partial_t u = \Delta u - g(u) \cdot \nabla u + f(u)   &\textrm{ in } \Omega \textrm{ for  } t>0, \\
			\sigma \partial_t u + \partial_\nu u =0  & \textrm{ on } \partial \Omega \textrm{ for  } t>0, \\
			u(\cdot,0) = u_0 \geq  0  &  \textrm{ in } \overline{\Omega},
    \end{array}
\right.		
\end{eqnarray}
where $g:\mathbb{R} \mapsto \mathbb{R}^N$, $f:\mathbb{R} \mapsto \mathbb{R}$, $\Omega$ is a bounded domain of $\mathbb{R^N}$ with $\mathcal{C}^2$-boundary $\partial \Omega$. We denote by $\nu:\partial \Omega \mapsto \mathbb{R}^N$ the outer unit normal vector field, and by $\partial_\nu $ the outer normal derivative.\\
\indent These equations arise in different areas, especially in population growth, chemical reactions and heat conduction. For instance, in the case of a heat transfer in a medium $\Omega$, the first equation $\partial_t u = \Delta u - g(u) \cdot \nabla u + f(u)$ is a heat equation including a nonlinear convection term $g(u)\cdot \nabla u$ and a nonlinear source $f$. On the boundary $\partial \Omega$, if $\sigma$ is positive, the dynamical boundary conditions describe the fact that a heat wave with the propagation speed $\frac{1}{\sigma}$ is sent into the region into an infinitesimal layer near the boundary due to the heat flux across the boundary (see \cite{CE} and \cite{Goldstein}).\\
\indent There are various results in the literature about the theory of blow-up for semilinear parabolic equations, in particular for reaction-diffusion equations, see e.g. \cite{Vazquez}, \cite{F}, \cite{FML}, and \cite{K}. In this work, we discuss a problem involving a nonlinear convection term. Whereas a Burgers' equation has been studied in \cite{BMR} in the one-dimensional case, we now consider a more general convection term and we set in a regular domain of $\mathbb{R}^N$.
After recalling some qualitative properties in Section \ref{S1}, we construct a global upper-solution for Problem \eqref{CRD} in Section \ref{SGE} and we deduce some conditions on $f$ and $g$ guaranteeing global existence of the solutions (Theorem \ref{GE_up}). In Section \ref{SBU}, we investigate two methods to ensure the blow-up of solutions of Problem \eqref{CRD}. The first one is an eigenfunction method valid for the model problem
\begin{eqnarray}\label{pbm2}
\left\{
    \begin{array}{ll}
			\partial_t u = \Delta u - g(u) \cdot \nabla u + u^p   &\textrm{ in } \overline{\Omega} \textrm{ for  } t>0, \\
			\sigma \partial_t u + \partial_\nu u =0  & \textrm{ on } \partial \Omega \textrm{ for  } t>0, \\
			u(\cdot,0) = u_0   &  \textrm{ in } \overline{\Omega},
    \end{array}
\right.		
\end{eqnarray}
with $p>1$ (Theorem \ref{BU_up}). We also derive some upper bounds for the blow-up time. The second method, devoted to the following problem
\begin{eqnarray}\label{pbm3}
\left\{
    \begin{array}{ll}
			\partial_t u = \Delta u - g(u) \cdot \nabla u + e^{pu}   &\textrm{ in } \overline{\Omega} \textrm{ for  } t>0, \\
			\sigma \partial_t u + \partial_\nu u =0  & \textrm{ on } \partial \Omega \textrm{ for  } t>0, \\
			u(\cdot,0) = u_0   &  \textrm{ in } \overline{\Omega},
    \end{array}
\right.		
\end{eqnarray}
with $p>0$, requires a self-similar lower-solution which blows up in finite time (Theorem \ref{BUthm3}). We prove the blow-up of solutions of  Finally, in Section \ref{SGO}, we determine the blow-up rate of the solutions of Problem \eqref{pbm2} in the $L^\infty$-norm when approaching the blow-up time (Theorem \ref{thmgo1}).\\
\indent Throughout, we shall assume the dissipativity condition
\begin{equation}\label{sigma0}
\sigma \geq 0 \textrm{ on } \partial \Omega  \times (0,\infty).
\end{equation}
In order to deal with classical solutions, we always assume that the parameters in the equations of Problem \eqref{CRD} are smooth
\begin{equation}\label{sigma1}
\sigma \in \mathcal{C}^1_b(\partial \Omega  \times (0,\infty) ) ,
\end{equation}
\begin{equation}\label{reactionf}
f \in \mathcal{C}^1(\mathbb{R} ) \ , \ f(s) >0 \textrm{ for } s> 0 ,
\end{equation}
and
\begin{equation}\label{convectiong}
g \in \mathcal{C}^1(\mathbb{R},\mathbb{R}^N ) .
\end{equation}
The initial data is continuous, non-trivial and non-negative in $\overline{\Omega}$
\begin{equation}\label{idata}
u_0 \in \mathcal{C}(\overline{\Omega}) ,  \ u_0 \not \equiv 0, \ u_0 \geq 0  .
\end{equation}
\indent Let $T = T(\sigma, u_0)$ denote the maximal existence time of the unique maximal classical solution of Problem \eqref{CRD}
\begin{displaymath}
u_\sigma \in \mathcal{C}(\overline{\Omega} \times [0,T) ) \cap \mathcal{C}^{2,1}(\overline{\Omega} \times (0,T) )
\end{displaymath}
with the coefficient $\sigma$ in the boundary conditions and the initial data $u_0$. As for the well-posedness and the local existence of the solutions of Problem \eqref{CRD}, we refer to \cite{vBDC}, \cite{CE} and \cite{E}. 
>From \cite{CE}, since the convection term depends linearly on the gradius $\nabla u$ of the solution, the maximal existence  time $T$ is the blow-up time of the solution with respect to the $L^\infty$-norm:
\begin{displaymath}
T = \inf \Big\{ \ s>0 \ \Big| \ \lim_{t \nearrow s}  \sup_{\overline{\Omega}} |u(x,t)|  = \infty \ \Big\}\ .
\end{displaymath}

\section{Qualitative properties}\label{S1}

The aim of this section is to compare the solutions for different parameters $\sigma$ and initial data $u_0$ and to summarize some positivity results on the classical solutions of Problem \eqref{CRD}. \\
\noindent Using the maximum principle from \cite{vBDC}, we extend some results obtained in \cite{vBPM2} in the case of reaction-diffusion to our problem with convection.

\begin{thm}\label{ut>d}
Assume hypotheses \eqref{sigma0} - \eqref{idata2}. Suppose that $\sigma$ does not depend on time 
\begin{equation}\label{dts0}
\sigma \in \mathcal{C}^1(\partial \Omega).
\end{equation}
 Then the solution $u$ of Problem \eqref{CRD} satisfies
\begin{displaymath}
u > 0  \ \textrm{ in } \ \overline{\Omega} \times (0,T(\sigma,u_0) ),
\end{displaymath}
\begin{displaymath}
\partial_t u \geq 0  \ \textrm{ in } \ \overline{\Omega} \times [0,T(\sigma,u_0) ),
\end{displaymath}
\begin{displaymath}
\partial_t u > 0  \ \textrm{ in } \ \overline{\Omega} \times (0,T(\sigma,u_0) ).
\end{displaymath}
Moreover, for all $\xi \in(0,T(\sigma,u_0))$, there exists $d>0$ such that
\begin{displaymath}
\partial_t u > d  \ \textrm{ in } \ \overline{\Omega} \times [\xi,T(\sigma,u_0) ).
\end{displaymath}
\end{thm}
\begin{pf}
Let $\tau \in (0,T(\sigma,u_0))$. Since $u$ is $\mathcal{C}^{2,1}(\overline{\Omega} \times [0,\tau])$ and because $f$ and $g$ are smooth (\eqref{reactionf} and \eqref{convectiong}), we can define these constants
\begin{displaymath}
C = \sup_{\overline{\Omega} \times [0,\tau]} g(u) \ \ \ \textrm{ and } \ \ \  M = \sup_{\overline{\Omega} \times [0,\tau]} g'(u) \cdot \nabla u - f'(u) .
\end{displaymath}
First, the positivity principle (Corollary 2.4 from \cite{vBDC}) applied to Problem \eqref{CRD} implies $u \geq 0$ in $\overline{\Omega} \times [0,\tau]$ since $f\geq 0$ by condition \eqref{reactionf}. Thus we obtain
\begin{eqnarray*}
\left\{
    \begin{array}{ll}
			\partial_t u \geq \Delta u - g(u) \cdot \nabla u \geq  \Delta u - C |\nabla u| &\textrm{ in } \Omega \textrm{ for  } t>0, \\
			\sigma \partial_t u + \partial_\nu u =0  & \textrm{ on } \partial \Omega \textrm{ for  } t>0, \\
			u(\cdot,0) = u_0  &  \textrm{ in } \overline{\Omega}.
    \end{array}
\right.		
\end{eqnarray*}
The strong maximum principle from \cite{vBDC} implies
\begin{displaymath}
m:=\min_{\overline{\Omega} \times [0,\tau]} u = \min_{\overline{\Omega} } u_0 \ ,
\end{displaymath}
and if this minimum $m$ is attained in $\overline{\Omega} \times (0,\tau]$, $u\equiv m$  in $\overline{\Omega} \times [0,\tau]$. Since $f>0$ in $(0,\infty)$, the first equation in Problem \eqref{CRD} leads to $m=0$, and we obtain $u_0\equiv 0$, a contradiction with equation \eqref{idata}. Hence $u>m\geq 0$ in $\overline{\Omega} \times (0,\tau]$. \\
Then, since the coefficients in the equations of Problem \eqref{CRD} are sufficiently smooth, classical regularity results in \cite{LSU} imply that $u \in \mathcal{C}^{2,2}(\overline{\Omega} \times [0,\tau])$. Thus $y=\partial_t u \in \mathcal{C}^{2,1}(\overline{\Omega} \times [0,\tau])$ and satisfies
\begin{eqnarray*}
\left\{
    \begin{array}{ll}
			\partial_t y = \Delta y - g(u) \cdot \nabla y - (g'(u) \cdot \nabla u )y + f'(u)y&\textrm{ in } \Omega \textrm{ for  } t>0, \\
			\sigma \partial_t y + \partial_\nu y =0  & \textrm{ on } \partial \Omega \textrm{ for  } t>0.
    \end{array}
\right.		
\end{eqnarray*}
By continuity, condition \eqref{idata2} implies $y(\cdot,0) \geq 0$ in $\overline{\Omega}$. Again, Corollary 2.4 from \cite{vBDC} implies $y\geq 0$ in $\overline{\Omega} \times [0,\tau]$. In order to apply properly the strong maximum principle, we have to introduce $w= y e^{Mt}\geq0$. By definition of $C$ and $M$, we obtain
\begin{eqnarray*}
\left\{
    \begin{array}{ll}
			\partial_t w \geq \Delta w - g(u) \cdot \nabla w  \geq \Delta w - C |\nabla w| &\textrm{ in } \Omega \textrm{ for  } t>0, \\
			\sigma \partial_t w + \partial_\nu w \geq 0  & \textrm{ on } \partial \Omega \textrm{ for  } t>0.
    \end{array}
\right.		
\end{eqnarray*}
Again, the strong maximum principle from \cite{vBDC} implies
\begin{displaymath}
\tilde{m}:=\min_{\overline{\Omega} \times [0,\tau]} w = \min_{\overline{\Omega} } w(\cdot,0) \ ,
\end{displaymath}
and if this minimum $\tilde{m}$ is attained in $\overline{\Omega} \times (0,\tau]$, $w\equiv \tilde{m}$  in $\overline{\Omega} \times [0,\tau]$. In particular, if $\tilde{m}=0$, we have $\partial_t u \equiv 0$ in $\overline{\Omega} \times [0,\tau]$, thus $u(\cdot,t)=u_0$ for all $t \in [0,\tau]$. Hence $u$ attains its minimum in $\overline{\Omega} \times (0,\tau]$, which is impossible according to the first part of the proof. Thus $w$ and $\partial_t u$ are positive in $\overline{\Omega} \times (0,\tau]$.\\
Finally, let $\xi \in (0,\tau)$. Because $y$ is continuous and thanks to the previous point, there exists $d>0$ such that $y(\cdot,\xi)>d$ in $\overline{\Omega}$. As $y$ satisfies
\begin{eqnarray*}
\left\{
    \begin{array}{ll}
			\partial_t y = \Delta y - g(u) \cdot \nabla y - \Big( g'(u) \cdot \nabla u + f'(u) \Big) y&\textrm{ in } \Omega \times [\xi,\tau], \\
			\sigma \partial_t y + \partial_\nu y =0  & \textrm{ on } \partial \Omega \times [\xi,\tau],
    \end{array}
\right.		
\end{eqnarray*}
the weak maximum principle from \cite{vBDC} implies 
\begin{displaymath}
\min_{\overline{\Omega} \times [\xi,\tau]} y = \min_{\overline{\Omega} } y(\cdot,\xi) \ .
\end{displaymath}
Hence $y>d$ in $\overline{\Omega}  \times [\xi,\tau]$. Note that $d$ depends only on $\xi$, not on $\tau$. Without this step, we only have $y \geq \tilde{m} e^{-M\tau}$ which may vanish as $\tau \to T(\sigma,u_0)$.\qed
\end{pf}

\noindent Let $0 \leq \sigma_1 \leq \sigma_2$ be two coefficients satisfying condition \eqref{sigma1}, $v_0 \leq u_0$ be two initial data fulfilling hypothesis \eqref{idata} and $w_0$ a function in $\mathcal{C}_0(\overline{\Omega})$ with $0 \leq w_0 \leq v_0$. Denote by $u_{\sigma_1}$, $u_{\sigma_2}$, $v$ and $w$ the maximal solutions of the following problems
\begin{eqnarray*}
\left\{
    \begin{array}{ll}
			\partial_t u_{\sigma_1} = \Delta u_{\sigma_1} - g(u_{\sigma_1}) \cdot \nabla u_{\sigma_1} + f(u_{\sigma_1})   &\textrm{ in } \Omega \textrm{ for  } t>0, \\
			\sigma_1 \partial_t u_{\sigma_1} + \partial_\nu u_{\sigma_1} =0  & \textrm{ on } \partial \Omega \textrm{ for  } t>0, \\
			u_{\sigma_1}(\cdot,0) = u_0  &  \textrm{ in } \overline{\Omega},
    \end{array}
\right.		
\end{eqnarray*}
\begin{eqnarray*}
\left\{
    \begin{array}{ll}
			\partial_t u_{\sigma_2} = \Delta u_{\sigma_2} - g(u_{\sigma_2}) \cdot \nabla u_{\sigma_2} + f(u_{\sigma_2})   &\textrm{ in } \Omega \textrm{ for  } t>0, \\
			\sigma_2 \partial_t u_{\sigma_2} + \partial_\nu u_{\sigma_2} =0  & \textrm{ on } \partial \Omega \textrm{ for  } t>0, \\
		u_{\sigma_2}(\cdot,0) = u_0   &  \textrm{ in } \overline{\Omega},
    \end{array}
\right.		
\end{eqnarray*}
\begin{eqnarray*}
\left\{
    \begin{array}{ll}
			\partial_t v = \Delta v - g(v) \cdot \nabla v + f(v)   &\textrm{ in } \Omega \textrm{ for  } t>0, \\
			\sigma_2 \partial_t v + \partial_\nu v =0  & \textrm{ on } \partial \Omega \textrm{ for  } t>0, \\
			v(\cdot,0) = v_0   &  \textrm{ in } \overline{\Omega},
    \end{array}
\right.		
\end{eqnarray*}
and
\begin{eqnarray*}
\left\{
    \begin{array}{ll}
			\partial_t w = \Delta w - g(w) \cdot \nabla w + f(w)   &\textrm{ in } \Omega \textrm{ for  } t>0, \\
			w =0  & \textrm{ on } \partial \Omega \textrm{ for  } t>0, \\
			w(\cdot,0) = w_0   &  \textrm{ in } \overline{\Omega}.
    \end{array}
\right.		
\end{eqnarray*}
Let $T(\sigma_1,u_0)$, $T(\sigma_2,u_0)$, $T(\sigma_2,v_0)$ and $T(w_0)$ be their respective maximal existence times.
For the reader convenience, we recall some results stemming from the comparison principle \cite{vBDC}.

\begin{thm}[\cite{vBPM}]\label{compare_thm}
Under the aforementioned hypotheses, we have
\begin{displaymath}
T(\sigma_2,u_0) \leq  T(\sigma_2,v_0)  \leq  T(w_0)
\end{displaymath}
and
\begin{displaymath}
0\leq w \leq v \leq u_{\sigma_2} \ \textrm{ in } \ \overline{\Omega} \times [0,T(\sigma_2,u_0) ) \ .
\end{displaymath}
In addition, if $u_0 \in \mathcal{C}^2(\overline{\Omega})$ with 
\begin{equation}\label{idata2}
\Delta u_0 - g(u_0)\cdot \nabla u_0 + f(u_0) \geq 0  \textrm{ in } \Omega ,
\end{equation}
we have
\begin{displaymath}
T(\sigma_1,u_0) \leq  T(\sigma_2,u_0)
\end{displaymath}
and
\begin{displaymath}
u_{\sigma_2} \leq u_{\sigma_1} \textrm{ in } \ \overline{\Omega} \times [0,T(\sigma_1,u_0) ) \ .
\end{displaymath}
\end{thm}

\noindent An important fact comes from the last statement of Theorem \ref{ut>d}. For any positive solution $u$ of Problem \eqref{CRD}, the maximum principle implies that for any $s \in (0 , T(\sigma,u_0)) $, there exists $c>0$  such that $u (\cdot,s) \geq c$ in $\overline{\Omega}$. Then, consider the solution $\tilde{u}$ of \eqref{CRD} with the constant initial data $c$ and $\tilde{\sigma} = \sup \sigma$ in the boundary conditions. Theorem \ref{compare_thm} implies $\tilde{u} \leq u$. Since $c$ satisfies equation \eqref{idata2}, Theorem \ref{ut>d} leads to $\partial_t \tilde{u} > d>0$. Thus, $\tilde{u}$ can be big enough after a long time (maybe it blows up). So does $u$, even if $u_0$ does not satisfy condition \eqref{idata2}.

\section{Global existence}\label{SGE}

In this section, we give some conditions on the function $g$ in the convection term, which ensure global existence of the solutions of Problem \eqref{CRD} for various reaction terms $f$. We use the comparison method from \cite{vBDC}. Thus, we just need to find an appropriate upper-solution of Problem \eqref{CRD} which does not blow up. This is our first lemma.

\begin{lem}\label{supersol}
Let $\alpha>0$ and $K>0$ be two real numbers and let $\eta \in \mathcal{C}^1([0,\infty))$ with $\eta' \geq \alpha^2$. For any integer $1 \leq j \leq N$, the function $U$ defined in $\Omega \times [0,\infty)$ by
\begin{displaymath}
U(x,t) = K \exp \Big( \alpha x_j + \eta(t) \Big) ,
\end{displaymath}
satisfies
\begin{eqnarray*}
\left\{
    \begin{array}{ll}
			\partial_t U \geq \Delta U - g(U) \cdot \nabla U + f(U)   &\textrm{ in } \Omega \textrm{ for  } t>0, \\
			\sigma \partial_t U + \partial_\nu U \geq 0  & \textrm{ on } \partial \Omega \textrm{ for  } t>0, \\
			U(\cdot,0) >  0  &  \textrm{ in } \overline{\Omega},
    \end{array}
\right.		
\end{eqnarray*}
if 
\begin{equation}\label{Conv>RD}
\alpha g_j (\omega) \geq \frac{f(\omega)}{\omega} \textrm{ for all  } \omega \geq 0
\end{equation}
and if 
\begin{equation}\label{sigmamin}
\sigma (x,t) \geq \frac{\alpha}{\eta'(t)} \textrm{ for all  } t>0.
\end{equation}
\end{lem}
\begin{pf}
A simple computation of the derivatives of $U$ leads us to 
\begin{displaymath}
\partial_t U - \Delta U + g(U) \cdot \nabla U  = \Big( \eta' - \alpha^2 \Big) U +  \alpha g_j(U) U \textrm{ in } \Omega \textrm{ for } t>0. 
\end{displaymath}
Since we assume $\eta' \geq \alpha^2$, hypothesis \eqref{Conv>RD} implies 
\begin{displaymath}
\partial_t U - \Delta U + g(U) \cdot \nabla U -f(U) \geq 0 \textrm{ in }  \Omega \times (0,\infty).
\end{displaymath}
Furthermore, on the boundary $\partial \Omega$, for $t>0$, we have
\begin{eqnarray}\label{U_bOm}
\sigma \partial_t U +\partial_\nu U & = &\Big( \sigma \eta'(t)+ \alpha \nu_j(x) \Big) U  \\
 &\geq &\Big( \sigma \eta'(t)- \alpha  \Big) U \geq 0 , \nonumber
\end{eqnarray}
by hypothesis \eqref{sigmamin} since $\nu$ is normalized, and clearly $U(x,0) =K \exp \Big( \alpha x_j + \eta(0) \Big) >0$ in $\overline{\Omega}$.\qed
\end{pf}

\begin{rmk}
In the case of the Dirichlet boundary conditions, we can use this upper-solution with the special choice $\eta \equiv 0$ (see \cite{QS}). But for the dynamical boundary conditions, we must use a positive time-dependent $\eta$ because our solutions are not bounded, see Theorem \ref{ut>d}. 
\end{rmk}

\noindent Now we can state the following theorems for a nonlinear reaction term $f$ growing as a power of $u$ (Problem \eqref{pbm2}), or as an exponential function (Problem \eqref{pbm3}).

\begin{thm}\label{GE_up}
Let $\sigma$ be a coefficient fulfilling conditions \eqref{sigma0},  \eqref{sigma1} and such that there exists $\delta >0$ with
$$ \inf_{\partial \Omega} \sigma   \geq  \delta \sup_{\partial \Omega} \sigma \textrm{ for } t>0
\ \textrm{ and }  \
\Big(\sup_{x\in \partial \Omega} \sigma(x,\cdot)\Big)^{-1} \in L^1_{\textrm{{\scriptsize loc}}}(\mathbb{R}^+) .$$ 
Assume $u_0$ satisfies condition \eqref{idata}. If there exists an integer $1 \leq j \leq N$ such that 
\begin{equation}\label{Conv>RD2}
\liminf_{\omega \to \infty} \frac{g_j(\omega)}{\omega^{p-1}} > 0, 
\end{equation}
then the solution of Problem \eqref{pbm2} is a global solution.
\end{thm}
\begin{pf}
In view of Theorem \ref{ut>d} and \eqref{Conv>RD2}, we can suppose that $u_0$ is sufficiently big such that there exists $C>0$ with
\begin{displaymath}
g_j(u) \geq C u^{p-1} \textrm{ in } \Omega  \textrm{ for } t>0.
\end{displaymath}
For $\displaystyle \eta (t) = C\delta^{-1} \int_0^t \Big( \sup_{x\in \partial \Omega} \sigma(x,s)\Big)^{-1} \ ds + C^2 t$, we have $\eta' \geq C^2$ and Equation \eqref{sigmamin} is satisfied. 
Let $K$ be  a positive number such that
\begin{displaymath}
K \geq u_0 (x) e^{-C x_j -\eta(0)} \textrm{ for all } x \in \overline{\Omega}.
\end{displaymath}
Then by hypotheses \eqref{sigma1},  \eqref{idata} and \eqref{sigmamin}, the function $U$ defined in Lemma \ref{supersol} is an upper-solution of Problem \eqref{pbm2} since $U(\cdot,0) \geq u_0$ in $\overline{\Omega}$. Using the comparison principle from \cite{vBDC}, the unique solution $u$ of Problem \eqref{CRD} satisfies
$$
0 \leq u(x,t) \leq U(x,t) \textrm{ for all } x \in \overline{\Omega} \textrm{ and } t>0 \ ,
$$
thus $u$ does not blow up.\qed
\end{pf}

\noindent This theorem holds in particular for a nonlinearity $g$ in the form $g(u)=( \alpha_1 u^{q_1} , \dots,$ $\alpha_i u^{q_i}, \dots, \alpha_N u^{q_N} )$ with at least one integer $j$ such that $\alpha_j >0$ and $q_j \geq p-1$. A similar result can be derived for Problem \eqref{pbm3}:

\begin{thm}
Under the aforementioned assumptions, the solution of Problem \eqref{pbm3} is a global solution if the convection term $g(u)\cdot \nabla u$ has (at least) one component $g_j$ satisfying $g_j(u)= \alpha_j  e^{q_j u}$ with $\alpha_j >0$ and $q_j>p$.
\end{thm}
\begin{pf}
Thanks to  $q_j>p$, condition \eqref{Conv>RD} is fulfilled because $\alpha_j  e^{q_j u} \geq \alpha_j  e^{p u}/u$ for $u$ sufficiently big.\qed
\end{pf}

\begin{rmk}
Condition \eqref{Conv>RD2} is optimal for Problem \eqref{pbm2}, see Theorems \ref{GE_up} and \ref{BU_up}. But it can be improved in some special cases, for example, if the reaction term is $f(u)=u \ln u$. Lemma \ref{supersol} implies that all solutions of Problem \eqref{CRD} are global if one component $g_j$ of $g$ satisfies $g_j(u) \geq \alpha_j \ln u$. In fact, in that case, every positive solution of \eqref{CRD} is global, without any assumption on the convection term $g$, since $\int_c ^\infty  \frac{1}{f(y)} \ dy = \infty$ for $c>0$, see Theorem 3.2 from \cite{CE}.
\end{rmk}

\noindent Condition \eqref{sigmamin} on $\sigma$ allows us to consider fast decaying functions $\sigma$, but, to ensure global existence, it is  essential that $\sigma$  does not vanish on the whole $\partial \Omega$. Indeed let us prove the following blow-up result related to the Neumann boundary conditions, for $\sigma \equiv 0$ on $\partial \Omega$.

\begin{thm}
Assume that $\sigma \equiv 0$, $u_0$ fulfills hypothesis \eqref{idata} and $f $ is positive in $(0,\infty)$ such that
\begin{equation}\label{NeumanBU}
\int_c^\infty \frac{1}{f(y)} \ dy < \infty \textrm{ for some } c>0 .
\end{equation}
Then every positive solution of Problem \eqref{CRD} blows up in finite time.
\end{thm}
\begin{pf}
Let $u$ be a non-trivial positive solution of
\begin{eqnarray}\label{NeuCRD}
\left\{
    \begin{array}{ll}
			\partial_t u = \Delta u - g(u) \cdot \nabla u + f(u)   &\textrm{ in } \Omega \textrm{ for  } t>0, \\
			 \partial_\nu u =0  & \textrm{ on } \partial \Omega \textrm{ for  } t>0, \\
			u(\cdot,0) = u_0  &  \textrm{ in } \overline{\Omega}.
    \end{array}
\right.		
\end{eqnarray}
Using the maximum principle from \cite{vBDC}, we have $u(\cdot,\xi)>0$ in $\overline{\Omega}$ for $\xi>0$. Hence, without loss of generality, we suppose $u_0 > c$ in $\overline{\Omega}$. Now, consider the maximal solution $z$ of the ODE $\dot{z} = f(z)$ with the initial data $\displaystyle z(0)= \inf \{ u_0(x) \ / \ x\in \overline{\Omega} \} $. Condition \eqref{NeumanBU} implies that its maximal existence time $T_z$ is finite:
\begin{displaymath}
T_z=\int_{z(0)} ^\infty \frac{1}{f(y)} \ dy < \infty .
\end{displaymath}
Since $\nabla z= 0$, $z$ is a lower solution of Problem \eqref{NeuCRD}. Using the comparison principle from \cite{vBDC}, we obtain $z(t) \leq u(\cdot,t)$ in $\overline{\Omega}$ for $t>0$. Thus, $u$ must blow up in finite time with $0<T<T_z$.\qed
\end{pf}

\begin{rmk}
This section illustrates the damping effect of the dissipative dynamical boundary conditions: we have shown that for nontrivial $\sigma \geq 0$ the maximal existence time of the solutions of Problem \eqref{CRD} can be strictly greater than the ones under the Neumann boundary conditions.
\end{rmk}

\section{Blow-up}\label{SBU} 

In this section, we investigate the blow-up in finite time for the solutions of Problems \eqref{pbm2} and \eqref{pbm3}. 
Let $G$ be a primitive of $g$ and suppose that there exist $\alpha >0$ and $q<p$ such that 
\begin{equation}\label{G<uq}
G(\omega) \leq \alpha \omega ^q \textrm{ for } \omega >0.
\end{equation}
By applying the eigenfunction method (see \cite{vBPM}, \cite{F} and \cite{K}), we obtain some conditions on the initial data $u_0$ which guarantee the finite time blow-up and we derive some upper bounds for the blow-up times. This is a general technique which can be applied to the following problem, where the boundary behaviour of the solutions is not involved:
\begin{eqnarray}\label{pbm2positive}
\left\{
    \begin{array}{ll}
			\partial_t u = \Delta u - g(u) \cdot \nabla u + u^p   &\textrm{ in } \overline{\Omega} \textrm{ for  } t>0, \\
			 u \geq 0  & \textrm{ on } \partial \Omega \textrm{ for  } t>0, \\
			u(\cdot,0) = u_0   &  \textrm{ in } \overline{\Omega}.
    \end{array}
\right.		
\end{eqnarray}
Henceforth, we denote by $\lambda$ the first eigenvalue of $-\Delta$ in $H_0^1(\Omega)$ and by $\varphi$ an eigenfunction associated to $\lambda$ satisfying
\begin{equation}\label{phi}
\varphi \in H_0^1(\Omega) ,\ 0<\varphi \leq 1  \textrm{ in }  \Omega.   
\end{equation}

\begin{thm}\label{thmbu1}
Let $\alpha>0$, $1<q<p$, $m=p/(p-q)$ and suppose $G$ satisfies condition \eqref{G<uq}. Assume hypotheses \eqref{sigma0} - \eqref{idata} are fulfilled. If 
\begin{equation}\label{eqphiu0}
\int_\Omega u_0\varphi^m\ dx > (2|\Omega|^{p-1}C)^{\frac{1}{p}}
\end{equation}
with
\begin{eqnarray*}
C= (p-1)|\Omega|\Big(\frac{4\lambda}{p-q}\Big)^\frac{1}{p-1} + \Big(\frac{4q}{p-q}\Big)^\frac{q}{p-q}\alpha^m \int_\Omega |\nabla \varphi|^m \ dx \ ,
\end{eqnarray*}
then the maximal classical solutions $u$ of Problem \eqref{pbm2positive} blow up in finite time $T$ satisfying
\begin{equation}\label{Tborne}
T\leq \frac{2 \int_\Omega u_0\varphi^m\,dx}{(p-1)\Big(|\Omega|^{1-p}\Big(\int_\Omega u_0\varphi^m\,dx\Big)^p-2C \Big)}=:\tilde T.
\end{equation}
\end{thm}
\begin{pf}
Define 
\begin{displaymath}
M(t)=\int_\Omega u(x,t)\varphi(x)^m\,dx.
\end{displaymath}
Thus,
\begin{displaymath}
\dot{M}(t)=\int_\Omega \Delta u\varphi^m\,dx-\int_\Omega g(u)\cdot\nabla u\,\varphi^m\,dx +\int_\Omega u^p\varphi^m\,dx.
\end{displaymath}
First, we prove that  
\begin{equation}\label{eqdelta2}
\int_\Omega \Delta u\varphi^m\,dx\geq -m\lambda|\Omega|^{\frac{p-1}{p}}\Big(\int_\Omega u^p\varphi^m\,dx\Big)^{\frac{1}{p}}.
\end{equation}
Observe that the behaviours of $\varphi$ and $\partial_\nu\varphi$ on $\partial\Omega$ imply
\begin{equation}\label{phinu}
\int_{\partial \Omega} \partial_\nu u \varphi^m \ ds = 0\ \textrm{ and } \ \int_{\partial \Omega} u\partial_\nu(\varphi^m) \ ds \leq 0 ,
\end{equation}
since $u\geq 0$ on $\partial\Omega$ for $t>0$. As in \cite{QS}, Equation \eqref{phinu} and Green's formula yield
\begin{equation}\label{modif_coro}
\int_\Omega \Delta u\varphi^m\,dx\geq -m\lambda \int_\Omega  u\varphi^m\,dx.
\end{equation}
Since $\varphi \leq 1$, $\int_\Omega u\varphi^m\,dx\leq\int_\Omega u\varphi^{\frac{m}{p}}\,dx$ and by H\"older's inequality, \eqref{eqdelta2} holds.\\
\noindent Now, we show that
\begin{equation}\label{nabla}
-\int_\Omega g(u)\cdot\nabla u\,\varphi^m\,dx \geq -m\alpha \Big(\int_\Omega |\nabla\varphi|^m  \,dx\Big)^{\frac{1}{m}} \Big(\int_\Omega u^p\varphi^m\,dx\Big)^{\frac{q}{p}}.
\end{equation}
By Green's formula and by definition of $G$ and $\varphi$, we have 
\begin{eqnarray*}
-\int_\Omega g(u)\cdot\nabla u\,\varphi^m\,dx&=&-\int_\Omega \textrm{div}(G(u))\varphi^m\,dx = m\int_\Omega (G(u)\cdot  \nabla\varphi) \varphi^{m-1} \,dx \ .
\end{eqnarray*}
Equation \eqref{G<uq} and H\"older's inequality lead to
\begin{eqnarray*}
\Big|\int_\Omega (G(u)\cdot\nabla\varphi) \varphi^{m-1} \,dx\Big| & \leq & \alpha \int_\Omega u^q\varphi^{m-1} |\nabla\varphi| \,dx \\
& \leq & \alpha \Big(\int_\Omega |\nabla\varphi|^m  \,dx\Big)^{\frac{1}{m}} \Big(\int_\Omega u^p\varphi^\frac{(m-1)p}{q}\,dx\Big)^{\frac{q}{p}} \ ,
\end{eqnarray*} 
and \eqref{nabla} is satisfied.\\ 
\noindent Henceforth, introduce
\begin{displaymath}
C_1=m\lambda|\Omega|^{\frac{p-1}{p}}\ \ \textrm{and}\ \ C_2=m\alpha \Big(\int_\Omega |\nabla\varphi|^m  \,dx\Big)^{\frac{1}{m}}.
\end{displaymath}
Then we obtain
\begin{equation}\label{eq1423}
\dot{M}(t)\geq\int_\Omega u^p\,\varphi^m\,dx-C_1\Big(\int_\Omega u^p \varphi^m \,dx\Big)^\frac{1}{p}-C_2\Big(\int_\Omega u^p \varphi^m \,dx\Big)^\frac{q}{p}.
\end{equation}
Set 
\begin{displaymath}
\varepsilon_1=\displaystyle\frac{p^\frac{1}{p}}{4^\frac{1}{p}C_1}\  \textrm{and}\ 
\varepsilon_2=\displaystyle\frac{p^\frac{q}{p}}{(4q)^\frac{q}{p}C_2}.
\end{displaymath}
Recall Young's inequality: for $a>0$ and $ \varepsilon>0$, $\displaystyle a= \frac{\varepsilon a}{\varepsilon} \leq \frac{\varepsilon^r a^r}{r} + \frac{1}{s\varepsilon^s}$ for $r,s>1$ with $r^{-1} + s^{-1} =1$. It yields
\begin{displaymath}
C_1 \Big( \int_\Omega u^p\,\varphi^m\,dx\Big)^\frac{1}{p} \leq\frac{1}{4}\int_\Omega u^p\,\varphi^m\ dx +\underbrace{\frac{p-1}{p \varepsilon_1^\frac{p}{p-1}}}_{:=C_3},
\end{displaymath}
and in the same way we have
\begin{displaymath}
C_2 \Big( \int_\Omega u^p\,\varphi^m\,dx\Big)^{\frac{q}{p}}\leq \frac{1}{4} \int_\Omega u^p\,\varphi^m\,dx+C_4,
\end{displaymath}
with 
\begin{displaymath}
C_4=\frac{1}{m \varepsilon_2^m}.
\end{displaymath}
Then
\begin{displaymath}
\dot{M}(t)\geq\frac{1}{2}\int_\Omega u^p\,\varphi^m\,dx-C
\end{displaymath}
with $C=C_3+C_4>0$.
By \eqref{phi} and H\"older's inequality, we obtain that
\begin{displaymath}
\dot{M}(t)\geq\frac{1}{2}|\Omega|^{1-p}M^p-C.
\end{displaymath}
Since $M$ is increasing with respect to $t$, owing to \eqref{eqphiu0} we have 
\begin{displaymath}
\dot{M}(t)\geq\Big( \frac{1}{2}|\Omega|^{1-p} -CM(0)^{-p} \Big) M^p ,
\end{displaymath}
and we can conclude that $u$ can not exist globally. To derive an upper bound for the blow-up time, we integrate the previous differential inequality between 0 and $t>0$. We obtain 
\begin{displaymath}
M(t)\geq \Bigg(M(0)^{1-p}-(p-1)\Big( \frac{1}{2}|\Omega|^{1-p} -CM(0)^{-p} \Big) t \Bigg)^{\frac{-1}{p-1}}.
\end{displaymath}
Hence $M$ blows up before $\tilde{T}= M(0)^{1-p}(p-1)^{-1} \Big( \frac{1}{2}|\Omega|^{1-p} -CM(0)^{-p} \Big)^{-1}$, so does $u$. Thus, $T\leq \tilde T$.\qed
\end{pf}

\noindent We can note that Condition \eqref{eqphiu0} on the initial data is only necessary to derive an upper bound for the maximal existence time. Thanks to Theorem \ref{ut>d}, we obtain:

\begin{thm}\label{BU_up}
Let $q<p$ and suppose $G$ satisfies 
$$
\limsup_{\omega \to \infty} \frac{G(\omega)}{\omega^q} < \infty.
$$
Assume that $\sigma$ and $u_0$ satisfy conditions \eqref{sigma0}, \eqref{sigma1} and \eqref{idata}. All the positive solutions of Problem \eqref{pbm2} blow up in finite time.
\end{thm}
\begin{pf}
Let $u$ be a positive solution of Problem \eqref{pbm2}. Theorem \ref{ut>d} permits to ensure that there exist $t_0>0$ and $C>0$ such that $u(\cdot,t_0)$ is big enough to satisfy Equation \eqref{eqphiu0} and  $G(u) \leq C u^q$ in $\Omega$ for $t>t_0$. Thus applying Theorem \ref{thmbu1} to $v(x,t)=u(x,t+t_0)$, we prove that $v$ blows up in a finite time $T_v$ satisfying \eqref{Tborne}. Hence, $u$ blows up in a finite time $T_u=t_0 +T_v$.\qed
\end{pf}

\noindent Now, we prove the blow-up of positive solutions of Problem \eqref{pbm3}. 

\begin{thm}\label{BUthm3}
Assume $\sigma$ and $u_0$ satisfy conditions \eqref{sigma0} - \eqref{idata}. If 
\begin{displaymath}
\limsup_{\omega \to \infty} \frac{|g(\omega)|}{e^{q\omega}} < \infty,
\end{displaymath}
then all the positive solutions of Problem \eqref{pbm3} blow up in finite time.
\end{thm}
\begin{pf}
Let $u$ be a positive solution of Problem \eqref{pbm3} and define $v=e^{\gamma u}$ with $\gamma \in (q,p)$ and $\gamma>1/2$. As in the previous proof, we suppose that $u$ is sufficiently big such that for some $C>0$ 
\begin{equation}\label{HBU2}
|g(u)| \leq C e^{qu} \textrm{ in } \Omega \textrm{ for } t>0. 
\end{equation}
Computing the derivatives of $v$, we obtain
\begin{displaymath}
\partial_ t v = \Delta v - \frac{1}{v} |\nabla v|^2 - g(u)\cdot \nabla v + \gamma v^{\frac{p+\gamma}{\gamma}} \textrm{ in } \Omega \textrm{ for } t>0.
\end{displaymath}
Using condition \eqref{HBU2}, we obtain
\begin{displaymath}
\partial_ t v   \geq\Delta v - \frac{1}{ v} |\nabla v|^2 - C v^\frac{q}{\gamma}|\nabla v| + \gamma v^{\frac{p+\gamma}{\gamma}} \textrm{ in } \Omega \textrm{ for } t>0.
\end{displaymath}
Young's inequality 
\begin{displaymath}
C v^\frac{q}{\gamma}|\nabla v| \leq \frac{C^2}{2} |\nabla v | ^2 + \frac{1}{2}  v  ^\frac{2q}{\gamma},
\end{displaymath}
leads to
\begin{displaymath}
\partial_ t v   \geq \Delta v - \frac{2+C^2}{2} |\nabla v|^2  + \gamma v^{\frac{p+\gamma}{\gamma}} - \frac{1}{2}v^{\frac{2q}{\gamma}} \textrm{ in } \Omega \textrm{ for } t>0,
\end{displaymath}
since $v\geq 1$. Morevover, we have 
\begin{displaymath}
\gamma v^{\frac{p+\gamma}{\gamma}} - \frac{1}{2}v^{\frac{2q}{\gamma}} \geq (\gamma -\frac{1}{2}) v^\frac{p+\gamma}{\gamma}
\end{displaymath}
by definition of $\gamma$. Thus, we obtain 
\begin{eqnarray}\label{pbmSW}
\left\{
    \begin{array}{ll}
			\partial_t v \geq \Delta v - \mu |\nabla v|^2 + \kappa v^\frac{p+\gamma}{\gamma}   &\textrm{ in } \overline{\Omega} \textrm{ for  } t>0, \\
			  v \geq 0  & \textrm{ on } \partial \Omega \textrm{ for } t>0, \\
			  v(\cdot,0)  > 0 & \textrm{ in } \overline{\Omega} ,
    \end{array}
\right.		
\end{eqnarray}
with $\mu =(2+C^2)/2$ and $\kappa=\gamma-1/2$. Without loss of generality (see Theorem \ref{ut>d}), we can suppose that $v(\cdot,0) \geq V(\cdot,0)$ in $\overline{\Omega}$, where 
\begin{displaymath}
V(x,t) = (1-\varepsilon t)^\frac{-1}{p-1} W\Big( \frac{ |x| }{(1-\varepsilon t)^m}\Big),
\end{displaymath}
with $0<m< \min \{ \frac{1}{2}, \frac{p-q}{q(p-1)} \}$, $W(y)= 1+A/2 - y^2/(2A)$, $A>\frac{1}{m(p-1)}$ and $\varepsilon < \frac{2\kappa (p-1)}{2+A}$. According to Souplet \& Weissler \cite{SW}, $V$  is a blowing-up sub-solution for Problem \eqref{pbmSW}. By the comparison principle from \cite{vBDC}, $v \geq V$ and $u$ blows up in finite time.\qed
\end{pf}

\begin{rmk}
In this section, we point out the accelerating effect of the dynamical boundary conditions, in comparison with the Dirichlet boundary conditions. Indeed, we prove that, even if the initial data $u_0$ is small, the solutions of Problem \eqref{pbm2} blow up in finite time. But, if we replace the dynamical boundary conditions by the Dirichlet boundary conditions in the second equation of Problem \eqref{pbm2}, it is well known that the solutions are global and decay to $0$ if the initial data are small enough, see for instance references \cite{Straughan} and \cite{Weissler}.
\end{rmk}

\section{Growth Order}\label{SGO}

In this section, we are interested in the blow-up rate for Problem \eqref{pbm2} when approaching the blow-up time $T$. For the convection term, we assume that
\begin{equation}\label{eqg2} 
g(u)=(g_1(u),\cdots,g_n(u)) \ \textrm{with}\ g_i(u)=u^q \ \forall i=1,\cdots,n, \ 1<q\in\,\mathbb{R}.
\end{equation} 

\noindent First, we derive a lower blow-up estimate for $p>q+1$, valid for any non-negative initial data $u_0\in \mathcal{C}(\overline{\Omega})$.

\begin{lem}
Let $p>q+1$, and assume hypotheses \eqref{sigma0} - \eqref{idata}. Then the classical maximal solution $u$ of Problem \eqref{pbm2} satisfies
\begin{displaymath}
\|u(\cdot,t)\|_\infty\geq (p-1)^{\frac{-1}{p-1}}\left(T-t\right)^{\frac{-1}{p-1}}
\end{displaymath}
for $0<t<T$.
\end{lem}
\begin{pf}
Let $t\in\,[0,T)$. Denote by $\zeta\in\,\mathcal{C}^1((0,t_1))$ the maximal solution of the IVP
\begin{eqnarray*}
\left\{\begin{array}{llll}
\dot{\zeta}&=&\zeta^p & \textrm{ in } (0,t_1)\\ 
\zeta(0)&=&\|u(\cdot,t)\|_\infty &
\end{array}
\right.
\end{eqnarray*}
with $t_1=\displaystyle\frac{1}{p-1}\|u(\cdot,t)\|_\infty^{1-p}$. Introduce $v \in \, \mathcal{C}(\overline{\Omega}\times[0,T-t)) \cap\mathcal{C}^{2,1}(\overline{\Omega}\times(0,T-t))$ defined by $v(x,s)=u(x,s+t)$ for $x \in \overline{\Omega}$ and $s\in[0,T-t)$. Then $v$ is the maximal solution of the problem
\begin{eqnarray*}
\left\{\begin{array}{lll}
\partial_tv=\Delta v-g(v)\cdot\nabla v+v^p & \textrm{ in } \ \Omega\ \textrm{ for }\ 0<s<T-t,\\
\sigma \partial_t v + \partial_\nu v =0 & \textrm{ on }  \partial\Omega\ \textrm{ for }\ 0<s<T-t,\\
v(\cdot ,0)=u(\cdot,t) & \textrm{ in }  \overline{\Omega}.\\
\end{array}\right.
\end{eqnarray*}
The comparison principle from \cite{vBDC} implies that $t_1 \leq T-t$.\qed
\end{pf}

\noindent This result remains valid for Problem \eqref{CRD} as soon as blow-up occurs. We just need a positive function $f$ such that an explicit primitive of $\frac{1}{f}$ is known.\\
\noindent We improve the technique developed in Theorem 2.3 in \cite{BMR} for an one-dimensional Burgers' problem and inspired by Friedman \& McLeod \cite{FML} to prove that the growth order of the solution of Problem \eqref{pbm2} amounts to $-1/(p-1)$ for $p>2q+1>3$, when the time $t$ approaches the blow-up time $T$.

\begin{thm}\label{thmgo1}
Suppose conditions \eqref{sigma0}, \eqref{idata}, \eqref{dts0} and \eqref{eqg2} are fulfilled. For
\begin{equation}\label{p>p*}
p>2q+1\ , 
\end{equation} 
there exists a positive constant $C$ such that the classical maximal solution $u$ of Problem \eqref{pbm2}
satisfies
\begin{displaymath}
\|u(\cdot,t)\|_\infty \leq \frac{C}{(T - t)^{1/p-1}} \quad for\,\,t\in [0,T).
\end{displaymath}
\end{thm}
\begin{pf}
Let $\beta>1$ such that 
\begin{equation}\label{eqpq}
p(p-1)(p-2q-1) = \frac{Nq^2}{\beta} >0 ,
\end{equation}
and choose $M>1$ such that 
\begin{displaymath}
M \geq \frac{Nq}{2(2q+1)}\beta^{\frac{2q}{p-2q-1}}.
\end{displaymath}
First, for $\xi\in\,(0,T)$, we shall prove that there exists $\delta>0$ such that
\begin{displaymath}
\partial_tu\geq\delta e^{-Mt}(u^p+\beta u^{2q+1})
\end{displaymath}
in $\overline{\Omega}\times[\xi,T)$.
Introduce 
\begin{displaymath}
J=\partial_tu-\delta d(t)k(u)
\end{displaymath}
with $d(t)=e^{-Mt}$ and $k(u)=u^p+\beta u^{2q+1}$. Note that classical regularity results from \cite{LSU} yield $J\in\ \mathcal{C}^{2,1}\left(\overline{\Omega}\times [\xi,T)\right)$. We recall that Theorem \ref{ut>d} implies that there exists $c>0$ such that $\partial_t u\geq c>0$ in $\overline{\Omega}\times[\xi,T)$. Thus, we can choose $\delta>0$ sufficiently small such that 
\begin{displaymath}
J(\cdot,\xi)\geq 0\ \quad \textrm{in}\ \ \overline{\Omega}.
\end{displaymath}
$J$ fulfills the boundary condition
\begin{displaymath}
\sigma \partial_t J + \partial_\nu J =\partial_t ( \sigma \partial_t u + \partial_\nu u) -\delta dk^{'}(u)( \sigma \partial_t u + \partial_\nu u)
-\sigma \delta d^{'}k(u) =\sigma \delta Me^{-Mt}k(u)\geq 0.
\end{displaymath}
Furthermore, $J$ satisfies
\begin{displaymath}
\partial_t J-\Delta J+g(u)\cdot\nabla J-(pu^{p-1}- g'(u)\cdot\nabla u)J=\delta dH(u)\textrm{ in } \overline{\Omega}\times [\xi,T),
\end{displaymath}
where
\begin{displaymath}
H(u):=pu^{p-1}k(u)-k^{'}(u)u^p+k^{''}(u)|\nabla u|^2-\frac{d^{'}}{d}k(u)-k(u) g'(u)\cdot\nabla u.
\end{displaymath}
To prove that $H(u)\geq 0$, we shall show that
\begin{eqnarray}\label{eqH}
\begin{array}{ll}
q\sqrt{N}u^{q-1}|\nabla u|(u^p+ \beta u^{2q+1}) \leq & M(u^p+ \beta u^{2q+1})+  \beta(p-2q-1)u^{p+2q}   \\
																							 &+(p(p-1)u^{p-2}+2q(2q+1)\beta u^{2q-1})|\nabla u|^2.
																							 \end{array}
\end{eqnarray}
Inequality \eqref{eqH} is trivial in the case where $M\geq q\sqrt{N}u^{q-1}|\nabla u|$. Now, suppose that 
$M< q\sqrt{N}u^{q-1}|\nabla u|$. When $q\sqrt{N}u^{q+1}\leq 2q(2q+1)|\nabla u|$, we have
$q\sqrt{N}u^{q-1}u^p|\nabla u|\leq p(p-1)u^{p-2}|\nabla u|^2$ and $q\sqrt{N}u^{3q}|\nabla u|\leq 2q(2q+1)u^{2q-1}|\nabla u|^2$ since $p>3$ then \eqref{eqH} follows. In the case where $q\sqrt{N}u^{q+1}> 2q(2q+1)|\nabla u|$, since 
\begin{displaymath}
u>\left(\frac{2(2q+1)}{Nq}M\right)^{\frac{1}{2q}}\geq\beta ^{\frac{1}{p-2q-1}},
\end{displaymath}
we obtain 
\begin{equation}\label{eq2546}
u^p+ \beta u^{2q+1}\leq 2 u^p.
\end{equation}
Moreover, \eqref{eqpq} yields
\begin{eqnarray*}
2 \sqrt{N}\,qu^{q+1}|\nabla u| & = & 2\sqrt{\beta p(p-1)(p-2q-1)}\,u^{q+1}|\nabla u|\\ 
& \leq & \left(\sqrt{\beta (p-2q-1)}\,u^{q+1}-\sqrt{p(p-1)}\,|\nabla u|\right)^2\\
 & &  + 2\sqrt{\beta p(p-1)(p-2q-1)}\,u^{q+1}|\nabla u|\\
& \leq &\beta (p-2q-1)u^{2(q+1)}+p(p-1)|\nabla u|^2.
\end{eqnarray*}
Thus, multiplying by $u^{p-2}$, we are led to 
\begin{displaymath}
2 \sqrt{N}qu^{q-1}|\nabla u|u^p\leq \beta(p-2q-1)u^{p+2q}+p(p-1)u^{p-2}|\nabla u|^2
\end{displaymath}
and by \eqref{eq2546}, the inequality \eqref{eqH} holds. Finally, we can conclude by the comparison principle from \cite{vBDC} that $J\geq 0$ in $\overline{\Omega}\times [\xi,T)$, in particular, $\partial_t u \geq \varepsilon u^p$ with $\varepsilon>0$.\\
Now, we shall derive the upper blow-up rate estimate of $\|u(\cdot,t)\|_\infty$ for $t\in\,[\xi,T)$. 
For each $x \in \Omega$, the integral
\begin{displaymath}
\int_t^\tau \frac{\partial_t u(x,s)}{u^p(x,s) }\,ds= \int_{u(x,t)}^{u(x,\tau)} \frac{1}{\eta^p}\,d\eta
\end{displaymath}
converges as $\tau \to T$. Integrating the inequality $\partial_t u \geq \varepsilon u^p$ leads to
\begin{displaymath}
\varepsilon (\tau -t) \leq \frac{u(x,\tau)^{1-p} - u(x,t)^{1-p}}{1-p} \leq \frac{u(x,t)^{1-p}}{p-1} \ .
\end{displaymath}
Letting $\tau \to T$ implies $u(x,t) \leq \Big( \varepsilon (p-1)(T-t) \Big)^\frac{-1}{p-1}$ and we can conclude as in the proof of Theorem 2.3 from \cite{BMR}.
\qed
\end{pf}

\section*{Acknowledgment}

The authors would like to thank Dr Mabel Cuesta for helpful discussions and valuable advices.

\bigbreak

\noindent{\sc Ga\"elle Pincet Mailly}\\
LMPA Joseph Liouville FR 2956 CNRS, Universit\'e Lille Nord de France\\ 50 rue F. Buisson, B.P. 699, F-62228 Calais Cedex, France \\ 
\it{e-mail: mailly@lmpa.univ-littoral.fr}
\smallskip

\noindent{\sc Jean-Fran\c cois Rault} \\
LMPA Joseph Liouville FR 2956 CNRS, Universit\'e Lille Nord de France\\
 50 rue F. Buisson, B.P. 699, F-62228 Calais Cedex, France \\
\it{e-mail: jfrault@lmpa.univ-littoral.fr}
\end{document}